\documentclass[12pt,a4paper,intlimits,sumlimits,onesided]{amsart}

\ifdefined\SMART
\usepackage[paperwidth=12cm,paperheight=24cm,left=5mm,right=5mm,top=10mm,bottom=5mm]{geometry}
\else
\calclayout
\fi

\usepackage[usenames,dvipsnames,svgnames,table,rgb]{xcolor}
\usepackage{amsfonts,amsmath,amsthm,amssymb,fancyhdr}
\usepackage[allcolors=blue,allbordercolors=blue,pdfborderstyle={/S/U/W 1}]{hyperref}
\usepackage[ocgcolorlinks]{ocgx2} 
\usepackage[integrals]{wasysym}
\usepackage{parskip} 
\usepackage{tcolorbox} 
\usepackage[inline]{enumitem} 
\usepackage{bbm} 
\usepackage{multicol} 

\usepackage[
backend=biber,
style=alphabetic,
sorting=nyt
]{biblatex}
\let\OLDthebibliography\thebibliography
\renewcommand\thebibliography[1]{
  \OLDthebibliography{#1}
  \setlength{\parskip}{0pt}
  \setlength{\itemsep}{0pt}
}

\usepackage{txfonts,pxfonts,tikz} 
\usepackage{tgschola}
\usepackage[T1]{fontenc}

\newcommand{\beql}[1]{\begin{equation}\label{#1}}
	\newcommand{\eeq}{\end{equation}}

\newcommand{\Abs}[1]{{\left|{#1}\right|}}

\newcommand{\Norm}[1]{{\left\|{#1}\right\|}}

\newcommand{\Set}[1]{{\left\{{#1}\right\}}}

\newcommand{\one}{{\mathbbm 1}}

\newcommand{\inner}[2]{{\langle #1, #2 \rangle}}

\newcommand{\vol}{{\rm vol\,}}
\newcommand{\ft}[1]{\widehat{#1}}

\newcommand{\ZZ}{\ensuremath{\mathbb{Z}}}

\newcommand{\RR}{\ensuremath{\mathbb{R}}}

\newcommand{\TT}{\ensuremath{\mathbb{T}}}

\newtheorem{theorem}{Theorem}

\newtheorem{conjecture}{Conjecture}

\theoremstyle{definition}

\begin{document}

\title{Fuglede's conjecture on orthogonal bases of exponentials}

\author{Mihail N. Kolountzakis}

\address{\href{https://math.uoc.gr}{Department of Mathematics and Applied Mathematics}, University of Crete, Voutes Campus, 700 13 Heraklion, Greece,\newline and \newline \href{https://ics.forth.gr/}{Institute of Computer Science}, Foundation of Research and Technology Hellas, N. Plastira 100, Vassilika Vouton, 700 13, Heraklion, Greece}

\email{kolount@uoc.gr}

\begin{abstract}
We discuss part of Fuglede's original paper \cite{fuglede1974operators} in which he posed his famous conjecture on which bodies in Euclidean space admit an orthogonal basis of exponentials for their $L^2$ space.
\end{abstract}


\maketitle

\section{Introduction}

In the 1970s Bent Fuglede \cite{fuglede1974operators} made the following conjecture, working on which has shaped many mathematicians' career.

\begin{conjecture}\label{conj:fuglede}
Suppose $\Omega\subseteq\RR^d$ is a bounded set of positive measure. Then $\Omega$ tiles space by translations if and only if the space $L^2(\Omega)$ has an orthogonal basis which consists of exponentials.
\end{conjecture}

\noindent{\em Tiling by translations.}
For $\Omega$ to tile space by translations means that there is a countable set $T \subseteq \RR^d$ such that
\begin{equation}\label{tiling}
\sum_{t \in T} \one_\Omega(x-t) = 1, \text{ for almost every } x \in \RR^d.
\end{equation}
In other words, we can translate the domain around so that the translates are disjoint and cover everything (up to measure $0$). This definition makes perfect sense with any nonnegative function $f$ in place of $\one_\Omega$. We say that $f\ge 0$ tiles with the set of translates $T$ (sometimes $T$ is called a \textit{tiling complement} of $f$ or $\one_\Omega$) at level $\ell$ if
\begin{equation}\label{f-tiling}
\sum_{t \in T} f(x-t) = \ell, \text{ for almost every } x \in \RR^d.
\end{equation}

\noindent{\em Spectrality.}
The exponential functions on $\RR^d$ are complex functions of the form $x \to e^{2\pi i \lambda \cdot x}$, where $\lambda\cdot x$ denotes the ordinary inner product on $\RR^d$. The vector $\lambda$ is called the \textit{frequency} of the exponential. If the set of exponentials with frequencies in $\Lambda\subseteq\RR^d$
$$
E(\Lambda) = \Set{x \to e^{2\pi i \lambda\cdot x}: \lambda\in\Lambda}
$$
is an orthogonal basis for $L^2(\Omega)$ we call $\Omega$ \textit{spectral} and $\Lambda$ its \textit{spectrum} (spectra are not unique if they exist as we can always translate $\Lambda$ without changing this property).

No relation between the set of tiling translates $T$ and the frequencies of the spectrum $\Lambda$ is claimed in Conjecture \ref{conj:fuglede}. But in one important (and easy) case, already proved in \cite{fuglede1974operators}, the relation is very clear.
\begin{theorem}\label{th:lattice}
$\Omega$ tiles by translations with the lattice $L$ if and only if it is spectral with the dual lattice $L^*$ as a spectrum.
\end{theorem}
If $L=A\ZZ^d$ is a lattice ($A$ is a nonsingular $d\times d$ matrix) its \textit{dual lattice} is $L^* = A^{-\top}\ZZ^d$. 

\section{Basic results on the Fuglede Conjecture}

Conjecture \ref{conj:fuglede} has given rise to hundreds of papers. Its variations are plentiful. It is easy to make sense of the conjecture not only in Euclidean spaces but in any locally compact abelian (LCA) group, such as $\ZZ^d$, $\TT^d$ and finite abelian groups such as cyclic groups. Tiling can be stated immediately in this context and the role of exponentials is played by the continuous characters of the group, the elements of the \textit{dual group} \cite{rudin1962groups}.

Early in the 21st century the ``spectral $\implies$ tiling'' direction for $\RR^d$ was disproved by Tao \cite{tao2004fuglede} and shortly thereafter the ``tiling $\implies$ spectral'' direction for $\RR^d$ was disproved by the author and Matolcsi \cite{kolountzakis2006tiles}. Both directions are now known to be false in $\RR^d$ and $\ZZ^d$ for $d\ge 3$ \cite{tao2004fuglede,kolountzakis2006hadamard,kolountzakis2006tiles,farkas2006tiles,farkas2006onfuglede} as well as for several classes of finite groups. It is also known to be true for some classes of finite groups (see for example \cite{kiss2021fuglede,kiss2022fuglede}). It is worth noting that all counterexamples for $\RR^d$ were first constructed as counterexamples in some finite groups, then lifted to $\ZZ^d$ and finally to $\RR^d$. The connections of Conjecture \ref{conj:fuglede} to number theory (especially in dimension $d=1$ or on the cyclic groups), combinatorics and, naturally, Fourier analysis are very pronounced \cite{kolountzakis2024orthogonal}. We should also mention that the question of spectrality makes sense also for measures. If $\mu$ is a probability measure on $\RR^d$, for example, when does it have a spectrum, i.e.\ an orthogonal basis for $L^2(\mu)$ which consists of exponentials? This is a very active area of research for which, however, we do not have a corresponding conjecture, as it is not clear how to define translational tiling by a measure.

\section{Spectrality as tiling}

Fuglede came to his conjecture while studying a problem posed to him by Segal and he proved Theorem \ref{th:lattice} above. He also proved that a triangle in the plane is not spectral. Additionally he claimed that the disk in the plane is not spectral but did not provide a full proof of this until much later \cite{fuglede-ball}. It should be intuitively clear to the reader that neither the triangle nor the disk can tile the plane by translations (though the triangle \textit{can} tile if reflections about the origin are allowed in the group of transformations along with translations).

Before discussing his proofs we shall make the following important observation \cite{kolountzakis2000packing} which was not made by Fuglede himself explicitly, though he was certainly aware of it.

Let $\Omega \subseteq \RR^d$ have finite Lebesgue measure and assume that the set $\Lambda \subseteq \RR^d$ is such that the set of the corresponding exponentials $E(\Lambda)$ is orthogonal in $L^2(\Omega)$. If $f \in L^2(\Omega)$ then Bessel's inequality gives
$$
\sum_{\lambda \in \Lambda} \Abs{\inner{f}{e_\lambda}_{L^2(\Omega)}}^2 \le \vol(\Omega) \cdot \Norm{f}_2^2.
$$
If $E(\Lambda)$ is also complete in $L^2(\Omega)$ then Bessel's inequality becomes an equality, and if this equality holds for all $f$ then we have completeness of $E(\Lambda)$.

Assuming $E(\Lambda)$ is an orthogonal basis we apply this now to $f(x) = e^{2\pi i t\cdot x}$ for an arbitrary frequency $t \in \RR^d$ and we obtain, after very little calculation,
\begin{equation}\label{spectrum-as-tiling}
\sum_{\lambda\in\Lambda} \Abs{\ft{\one_\Omega}}^2 (t-\lambda) = \vol(\Omega)^2.
\end{equation}
Equation \eqref{spectrum-as-tiling} clearly represents a tiling situation as in \eqref{f-tiling} where the tile is now the nonnegative function $\Abs{\ft{\one_\Omega}}^2 \in L^1(\RR^d)$ and the translation set is the spectrum $\Lambda$. Since the trigonometric polynomials are dense in $L^2(\Omega)$ we can see that \eqref{spectrum-as-tiling} implies $\sum_{\lambda \in \Lambda} \Abs{\inner{f}{e_\lambda}}^2 = \vol(\Omega) \Norm{f}_2^2$ for any $f \in L^2(\Omega)$. Clearly Conjecture \ref{conj:fuglede} is then equivalent to the more symmetric form:
\begin{quotation}
$\Omega$ tiles by translations $\iff$ $\Abs{\ft{\one_\Omega}}^2$ tiles by translations.
\end{quotation}
Since $\inner{e^{2\pi i \lambda\cdot x}}{e^{2\pi i \mu\cdot x}}_{L^2(\Omega)} = \ft{\one_\Omega}(\lambda-\mu)$ we also have the important observation that orthogonality of $E(\Lambda)$ is equivalent to the difference set $\Lambda-\Lambda$ being a subset of the zero set $\Set{\ft{\one_\Omega}=0} \cup \Set{0}$ of the Fourier Transform of the indicator function of $\Omega$.

\section{The triangle is not spectral}

\begin{figure}[ht]
\ifdefined\SMART\resizebox{10cm}{!}{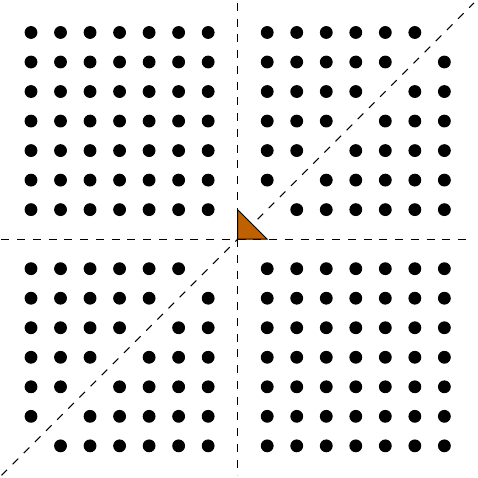}\else
\input zeros.pdf_t
\fi
\caption{The zero set $Z$ of the Fourier Transform of the triangle shown with one vertex at the origin and the others at $e_1, e_2$. The dots (zeros) are on integer points. The three dashed lines are missing (they are non-zeros).}\label{fig:zeros}
\end{figure}

Fuglede proved that the triangle $\Omega$ is not spectral by first finding the zeros of the Fourier Transform of $\one_\Omega$. With the Fourier Transform normalization
$$
\ft{f}(\xi) = \int_{\RR^d} f(x) e^{-2\pi i \xi\cdot x}\,dx
$$
and working on the specific triangle with vertices $(0, 0), (1, 0)$ and $(0, 1)$ (we can always apply a linear transformation to our domain as this does not change spectrality) Fuglede computed that the zero set of the Fourier Transform of the triangle is the set
$$
Z = \Set{(m, n)\in\ZZ^2: m \neq 0,\ n\neq 0,\ m\neq n},
$$
shown in Figure \ref{fig:zeros}.

At this point it is clear to us now (but, apparently, was not clear to Fuglede back then, judging  from his proof) that there is no spectrum. We now know from \eqref{spectrum-as-tiling} that any possible spectrum $\Lambda$ must have density (number of points in a large disk, divided by the area of the disk) equal to $\vol(\Omega)=\frac12$. This is true for any tiling: the volume of the tile (or the integral of the tiling function, as in \eqref{spectrum-as-tiling}) times the density of the translation set must equal the level of the tiling. Since we are free to translate any possible spectrum $\Lambda$ we assume, as we often do, that $0 \in \Lambda$. Orthogonality then dictates that $\Lambda \subseteq Z \subseteq \ZZ^2$ and the shape of $Z$ shows that there is at most one point of $\Lambda$ on any one horizontal (or vertical, or parallel to the line $y=x$) line. It follows that $\Lambda$ has $O(R)$ points in any ball of radius $R$ so it has zero density, not density $1/2$ as it should have.

On the other hand, as Fuglede \cite{fuglede1974operators} pointed out, there are infinite orthogonal sets $\Lambda$ for the triangle. For instance take $\Lambda = \Set{(-n, n):\ n\in\ZZ}$.

Let us now describe an alternative proof that the triangle is not spectral which avoids the precise calculation of the zero set of the Fourier Transform and which is much more general \cite{kolountzakis2002class}. Let us write $f(x, y)=\one_\Omega(x, y)$ for the indicator function of the triangle. To compute its Fourier Transform restricted on the $x$-axis we can, using Fubini's theorem, first project $f$ onto the $x$-axis and take the projection's one-dimensional Fourier Transform.

This projection $g(x) = \int_{\RR} f(x, y)\, dy$ is the piecewise linear function which is equal to 1 at 0, equal to 0 at 1 and is 0 off $[0, 1]$. To find the roots of $\ft g$ is the same as finding the roots of $\ft{g'}(\xi) = 2\pi i \xi \ft g(\xi)$ except at 0. The derivative of $g$ is the measure
$$
\mu = \delta_0 - \one_{[0, 1]},
$$
whose Fourier Transform is
$$
\ft\mu(\xi) = 1 - \ft{\one_{[0,1]}}(\xi).\ \ \ (\xi\in\RR)
$$
The second term is the Fourier Transform of an integrable function and so it tends to 0 as $\Abs{\xi}\to\infty$. This implies that for some constant $K>0$ we have $\Abs{\ft{\mu}(\xi)}>\frac12$ if $\Abs{\xi} \ge K$ and
\begin{equation}\label{ft-bound}
\Abs{\ft{\one_\Omega}(\xi, 0)} = \Abs{\ft{g(\xi)}} \ge \frac{c}{\Abs{\xi}} \text{ for } \Abs{\xi} \ge  K \text{ and some constant } c>0.
\end{equation}
It is easy to see using the divergence theorem that for $\zeta=(\xi, \eta) \in \RR^2\setminus\Set{0}$ we have
$$
\ft{\one_\Omega}(\zeta) = -{1\over i\Abs{\zeta}}
        \int_{\partial\Omega} e^{-2\pi i \zeta \cdot t}
         {\zeta\over\Abs{\zeta}}\cdot{\nu(t)}\,d\sigma(t),
\ \ \ x\neq 0,
$$
where $\nu(t) = (\nu_1(t),\nu_2(t))$ is the outward unit normal vector to $\partial\Omega$ at $t \in \partial\Omega$ and
$d\sigma$ is the arc-length measure on $\partial\Omega$. From this it follows that
\begin{equation}\label{grad-ft-bound}
\Abs{\nabla\ft{\one_\Omega}(\zeta)} \le \frac{C}{\Abs{\zeta}}, \text{ for } \Abs{\zeta}\ge 1.
\end{equation}
It then follows from \eqref{ft-bound} and \eqref{grad-ft-bound} that for some $\epsilon>0$ and for all $\xi, \eta$ with $\Abs{\xi}\ge K$ and $\Abs{\eta} < \epsilon$ the Fourier Transform $\ft{\one_\Omega}(\xi, \eta) \neq 0$.

\begin{figure}[ht]
\ifdefined\SMART\resizebox{10cm}{!}{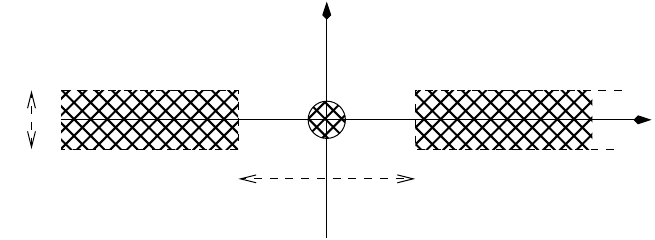}\else
\input no-zeros.pdf_t
\fi
\caption{The shaded region has no zeros of $\ft{\one_\Omega}$.} \label{fig:no-zeros}
\end{figure}

It is easy to see that for some positive $\rho$ (the magnitude of the smallest root of $\ft{\one_\Omega}$) we have $\Abs{\lambda-\mu} \ge \rho$ for any two distinct $\lambda, \mu \in \Lambda$. From the zero-free region (see Fig.\ \ref{fig:no-zeros}) we established above it follows easily that the set $\Lambda$ has at most a bounded number of points in each horizontal strip of width $2\epsilon$, which implies that $\Lambda$ has density 0 and cannot therefore be a spectrum.

\section{The disk is not spectral}

In \cite{fuglede1974operators} Fuglede also stated that the disk is not a spectral set. This, along with his result about the triangle, matched his expectation that is contained in Conjecture \ref{conj:fuglede}. It is clear that the disk cannot tile and it is well known since the time of Minkowski that every convex set that tiles by translation must be symmetric (this latter property was extended to spectral sets early on \cite{kolountzakis2000nonsymmetric}: spectral convex sets must be symmetric). The impossibility of translational tiling with a triangle should be intuitevely obvious: how are we to cover the outside of one of the triangle's edges, say edge $e$? There is no other edge of the triangle that is parallel to $e$ and facing in the opposite direction.

Fuglede outlined (in very broad strokes) how the proof that the disk is not spectral should proceed. In fact he claimed that using properties of the zeros of the Bessel function $J_1$ (essentially the Fourier Transform of the disk) one can prove that there is no infinite orthogonal set $\Lambda$ (which is clearly stronger than non-spectrality).

Fuglede returned later to this problem \cite{fuglede-ball} and proved that there are no infinite sets of orthogonal exponentials for a ball in any dimension. He used, as he predicted, asymptotic properties of roots of Bessel functions. The same result was proved also in \cite{iosevich-katz-pedersen,iosevich-rudnev-smooth-bodies}. See also \cite[Section 3.4]{kolountzakis2004} for an almost ``accidental'' proof that the disk in the plane is not spectral.

Let us make clear here that, although we know that any orthogonal set of exponentials for the disk must be finite, we do not know that the cardinality of all such orthogonal sets is uniformly bounded. (In other words, it is still conceivable that there are arbitrarily large orthogonal sets for the disk.) Fuglede himself expected the uniform bound to be 3 (which is clearly achievable for the disk by taking three points at an equilateral triangle whose side-length is a root of the Fourier Transform of the disk). Thus it makes sense to ask, for an orthogonal set $\Lambda$ for the disk, how many points it can have in a disk of radius $R$. The first result of this type was that this quantity is $O(R)$ \cite{iosevich-jaming}. This was improved to $O(R^{2/3})$ in \cite{iosevich2013size} and this upper bound stood until very recently when it was improved to $O(R^{3/5+\epsilon})$ for any $\epsilon>0$ \cite{zakharov2024sets}.

We would like to show here a method to show nonspectrality of the disk \cite{kolountzakis2004distance} that uses some Ramsey properties of sets of positive Lebesgue density. (It is however known by now that Conjecture \ref{conj:fuglede} holds for convex bodies \cite{lev2022fuglede} so we are not proving anything new.)

\begin{figure}[ht]
\ifdefined\SMART\resizebox{10cm}{!}{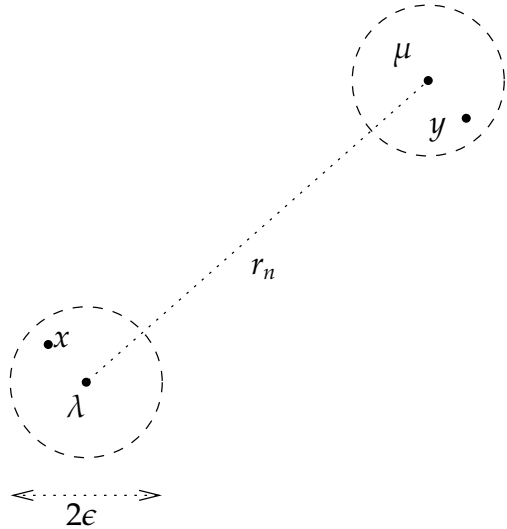}\else
\input distances.pdf_t
\fi
\caption{The set $E$ consists of an $\epsilon$-ball around each point of the spectrum.} \label{fig:distances}
\end{figure}

There are two main ingredients in this proof. Both hold for smooth symmetric convex bodies of positive curvature everywhere, and so does the resulting non-spectrality, but we will restrict the discussion to the ball to keep it simple.
\begin{enumerate}
\item
The zeros of $\ft{\one_\Omega}$, where $\Omega$ is the unit ball, are spheres, centered at the origin whose radii $r_n$ are asymptotic to an arithmetic progression
\begin{equation}\label{asymptotics}
r_n = A + Bn + o(1)
\end{equation}
where $A, B$ are constants (see e.g.\ \cite{iosevich2013size}).
\item 
If a set $E \subseteq \RR^d$ has positive density with respect to Lebesgue measure (this means that there is a constant $c>0$ and arbitrarily large balls such that in each of these balls the set $E$ takes up proportion at least $c$ of their volume) then all sufficiently large distances are realized within the set $E$ \cite{furstenberg1990ergodic,bourgain1986szemeredi,kolountzakis2004distance}.
\end{enumerate}
Suppose then that $\Lambda \subseteq \RR^d$ is a spectrum for the unit ball. Then any two points of $\Lambda$ are distance at least $C>0$ apart, where $C$ is a constant and $\Lambda$ has positive (counting) density. Let $\epsilon = \frac1{10}\min\Set{B, C}$ and define then the set
$$
E = \Lambda+B_{\epsilon},
$$
which has positive density wit respect to Lebesgue measure since $\Lambda$ has positive counting density and the balls $B_\epsilon$ of radius $\epsilon$  centered at $\Lambda$ are disjoint. Take two points $x, y \in E$ (see Fig.\ \ref{fig:distances}). Then there are two points $\lambda, \mu \in \Lambda$ such that
$$
\Abs{x-\lambda}<\epsilon \text{ and } \Abs{y - \mu} <\epsilon,
$$
which implies that
$$
\Abs{ \, \Abs{x-y} - \Abs{\lambda-\mu} \,} < 2\epsilon.
$$
But $\Abs{\lambda-\mu}$ is equal to some $r_n$ and so very close to a number of the form $A+Bn$ from \eqref{asymptotics} if only $\Abs{\lambda-\mu}$ is large enough. It follows that $\Abs{x-y}$ is within $2\epsilon$ from some $r_n$. Since the distance between successive $r_n$s is $B+o(1)$, which is much larger than $2\epsilon$, it follows that the possible values for $\Abs{x-y}$ have gaps going all the way to infinity, contradicting the fact, mentioned above, that all sufficiently large distances occur between points of $E$. We have arrived at a contradiction, so the ball does not have a spectrum (it does not even have an orthogonal set $\Lambda$ of positive counting density, by the argument we just presented).
 
\printbibliography

\end{document}